\documentclass[11pt,a4paper,reqno,oneside]{amsart}
\usepackage[noBBpl]{mathpazo}

\usepackage{config, cancel, xfrac, empheq}

\usepackage[most]{tcolorbox}
\usepackage{the}

\usepackage{nameref}
\newcounter{mylabelcounter}
\makeatletter
\newcommand{\labelText}[2]{%
#1\refstepcounter{mylabelcounter}%
\immediate\write\@auxout{%
  \string\newlabel{#2}{{1}{\thepage}{{\unexpanded{#1}}}{mylabelcounter.\number\value{mylabelcounter}}{}}%
}%
}
\makeatother

\makeatletter
\providecommand*{\twoheadrightarrowfill@}{%
  \arrowfill@\relbar\relbar\twoheadrightarrow
}
\providecommand*{\xtwoheadrightarrow}[2][]{%
  \ext@arrow 0579\twoheadrightarrowfill@{#1}{#2}%
}
\makeatother

\begin{document}
\begin{abstract}
We give a new elementary proof of the theorem that a natural map from Milnor's construction \(F[S^1]\) to the simplicial group \(\mathrm{AP}\) of pure braids is injective. Our approach is group-theoretic and does not rely on Lie algebras.
\end{abstract}

\title{On injectivity of Cohen--Wu homomorphism}
\author{Vasily Ionin}
\address{St. Petersburg Department of Steklov Institute of Mathematics, Fontanka 27, St. Petersburg 191011, Russia}
\email{ionin.code@gmail.com}

\maketitle

\tableofcontents

\section{Introduction}

For every pointed simplicial set \(X\), there is Milnor's free group construction \(F[X]\) whose geometric realization is homotopy equivalent to \(\Omega \Sigma |X|\). In particular, \(F[S^n]\) is a simplicial group such that \(F[S^n]_k\) is the free group of rank \(\binom{n}{k}\) and \(|F[S^n]| \simeq \Omega S^{n+1}\).

Let \(M\) be a smooth connected surface. Denote by \(P_n(M)\) the pure braid group with \(n\) strands on the surface \(M\). When \(M\) is a regular disk, we denote the pure braid group simply as \(P_n\). The sequence \(\P(M) \defeq \{P_{n+1}(M)\}_{n \ge 0}\) forms a \(\Delta\)-group with face maps given by string deletions.
If \(M\) admits a non-degenerate vector field, then \(\P(M)\) could be equipped with the degeneracies (given by string doublings) and become a simplicial group.
For instance, there is a contractible simplicial group \(\AP \defeq \P(\RR^2)\) such that \(\AP_n = P_{n+1}\). The 1-simplex \(A_{1,2} \in \AP_1 = P_2\) is spherical, which defines a simplicial map \(S^1 \to \AP\). There is a unique homomorphism of simplicial groups
\begin{equation}
\Theta \colon F[S^1] \to \AP
\end{equation}
that extends this map.
In \cite{CW04, CW10}, Fred Cohen and Jie Wu prove the following theorem.
\begin{Theorem}{Cohen-Wu embedding}{cohen-wu-embedding}
The homomorphism \(\Theta\) is injective.\\
Hence, the homotopy type of the coset simplicial set \(\AP / F[S^1]\) is that of the 2-sphere \(S^2\).
\end{Theorem}
This theorem is of great interest because it realizes homotopy groups \(\pi_n(S^2)\) as some natural sub-quotients of pure braid groups; in particular it allows to represent elements of homotopy groups of 2-sphere with braids.

The lower central series of a group \(G\) is defined inductively as
\[\gamma_1(G) \defeq G, \gamma_{n+1}(G) \defeq [\gamma_n(G), G].\]
With a group \(G\) one can associate a graded Lie ring \(\mathrm{gr}(G)\) given by
\[\mathrm{gr}(G) \defeq \bigoplus_{n \ge 1} \frac{\gamma_n(G)}{\gamma_{n+1}(G)},\]
with a Lie bracket given by a commutator \([x, y] \defeq xyx^{-1}y^{-1}\). If \(f \colon G \to H\) is a homomorphism of groups, and \(G\) is residually nilpotent, then the injectivity of \(f\) can be deduced from the injectivity of the associated map \(\mathrm{gr}(f) \colon \mathrm{gr}(G) \to \mathrm{gr}(H)\). In \cite{CW04, CW10}, authors prove that \(\Theta_n\colon F[S^1]_n \to P_{n+1}\) is injective by showing that the Lie algebra homomorphism \(\mathrm{gr}(\Theta_n)\) is injective. Their proof contains heavy computations that rely on the explicit presentation of the Lie algebra \(\mathrm{gr}(P_{n+1})\) given by `infinitesimal braid relations'. It feels confusing that one has to descend to the level of Lie algebras, since the pure braid group is replete with free subgroups which are free for obvious reasons. For example, one can show very easily that the subgroup \(\ker(d_i)\) in \(P_{n+1}\), consisting of braids that become trivial after removing the \(i\)th strand, is free of rank \(n\).

In this paper, we give a new proof of Theorem~\ref{t:cohen-wu-embedding}, which is internal to the group theory.

Our strategy is straightforward. We have to show that some particular braids generate a free subgroup in a pure braid group. First, we `untwist' the problem by picking up a convenient generating set and by finding an automorphism of pure braid group such that it maps generators of subgroup to especially simple-looking braids.
Second, we solve the resulting group-theoretic problem with standard techniques (HNN extensions).

\medskip
\textbf{Acknowledgment.}
I am deeply grateful to Vanya Vasilyev and Artem Semidetnov for hours of discussion and general support. I would like to thank Ivan Gaidai-Turlov, Ilya Alexeev, Lev Mukoseev, and Roman Mikhailov for the early listening of proof attempts and helpful insights.

This work was supported by the Ministry of Science and Higher Education of the Russian Federation (agreement \texttt{075-15-2025-344} dated 29/04/2025 for Saint Petersburg Leonhard Euler International Mathematical Institute at PDMI RAS). The work was supported by the Theoretical Physics and Mathematics Advancement Foundation ``BASIS'', project no. \texttt{24-7-1-26-3}.

\section{Preliminaries}
We refer to \cite{MK99} for standard notions of braid theory. Recall that the Artin's pure braid group \(P_n\) is generated by the braids \(A_{i,j}, 1 \le i < j \le n\). For example, \(P_2\) is an infinite cyclic group generated by \(A_{1,2}\).
Given \(1 \le i < j \le n\), let \(A_{j,i} \defeq A_{i,j}\), and let \(A_{i,i} \defeq 1\).

Given \(1 \le i \le n\), let
\[A_{0,i} \defeq (A_{1,i} A_{2,i} \ldots A_{n,i})^{-1}.\]
These braids are important since \(\langle A_{0,1}, A_{0,2}, \ldots, A_{0,n} \rangle\) is the group of spherically trivial braids (that is, the kernel of the natural homomorphism \(P_n \to P_n(S^2)\)).

For \(1 \le i \le n\), let \(t_i, \widetilde A_{0,i} \in P_n\) be defined by the following formulas:
\begin{align*}
t_i &\defeq A_{1,i} A_{2,i} \ldots, A_{i-1,i}, \\
\widetilde A_{0,i} &\defeq A_{0,i} \cdot t_i^2.
\end{align*}
The braids \(\widetilde A_{0,i}\) are important, as they turn out to form a generating set for \(\mathrm{Im}(\Theta_{n-1})\).

\begin{figure}[H]
\centering
\includegraphics[width=\textwidth]{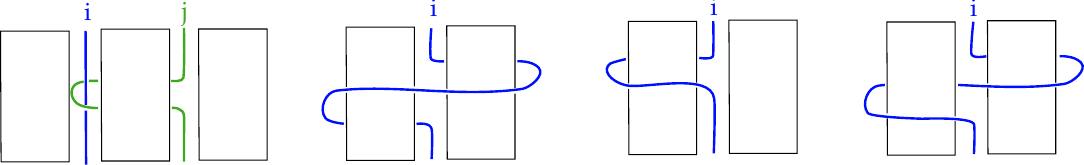}
\caption{Pictures of the geometric braids \(A_{i,j}, A_{0,i}, t_i\), and \(\widetilde A_{0,i}\).}
\end{figure}

\begin{Proposition}{}{bigprod}
One has
\begin{align*}
A_{0,1} A_{0,2} \ldots A_{0,n} &= z_n^{-2}, \\
\widetilde A_{0,1} \widetilde A_{0,2} \ldots \widetilde A_{0,n} &= 1.
\end{align*}
\end{Proposition}
\begin{proof}
The first equation can be found in \cite[equation~(25)]{AIM25}.

The second one can be proved in a similar way (one shows that this braid is in the center of the group \(P_n\), and notices that its image in the abelianization of \(P_n\) is zero).
\end{proof}

Let \(\mu_n\) be an automorphism of \(P_n\) given by the formula
\[\mu_n(A_{i,j}) = A_{i-1,j-1}, \quad 1 \le i < j \le n.\]

The following computation is straightforward.
\begin{Lemma}{}{mucomp}
One has
\[\mu_n(A_{0,i}) = \begin{cases}
A_{i-1,n}, &2 \le i \le n, \\
A_{0,n} z_n^2, & i = 1.
\end{cases}\]\qed
\end{Lemma}

For every \(1 \le k \le n\) there is a homomorphism of string deletion \(d_k \colon P_n \to P_{n-1}\)
given by
\[d_k(A_{i,j}) = \begin{cases}
A_{i, j}, & 1 \le i < j < k, \\
A_{i, j-1}, & 1 \le i < k < j \le n, \\
\triv, & k = i < j \le n \;\text{or}\; 1 \le i < j = k, \\
A_{i-1,j-1}, & k < i < j \le n. \\
\end{cases}\]
For every \(1 \le k \le n\) there is a homomorphism of string duplication \(s_k \colon P_n \to P_{n+1}\)
given by
\begin{align*}
s_k(A_{i,j}) &= \begin{cases}
A_{i, j}, & 1 \le i < j < k, \\
A_{i,j}A_{i,j+1}, & 1 \le i < j = k, \\
A_{i, j+1}, & 1 \le i < k < j \le n, \\
A_{i, j+1}A_{i+1,j+1}, & k = i < j \le n, \\
A_{i+1,j+1}, & k < i < j \le n.
\end{cases}
\end{align*}

A simplicial group structure on \(\AP\) is given by \(d^\AP_* \defeq d_{*+1}\) and \(s^\AP_* \defeq s_{*+1}\).

\section{Reduction to the group theoretic problem}
Recall that for every \(n \ge 1\) the group \(F[S^1]_n\) is free of rank \(n\) generated by simplices
\[z_i = \Big(\underbrace{s^{S^1}_0 \circ \ldots \circ s^{S^1}_0}_{i-1} \circ \underbrace{s^{S^1}_1 \circ \ldots \circ  s^{S^1}_1}_{n-i}\Big)(z), \quad 1 \le i \le n,\]
where \(z\) is the unique non-degenerate \(1\)-simplex of \(S^1\).
Since \(\Theta_n(z) = A_{1,2}\), the homomorphism \(\Theta_n \colon F(z_1, \ldots, z_n) \to P_{n+1}\) is given by
\[\Theta_n(z_i) = \left(\underbrace{s^\AP_0 \circ \ldots s^\AP_0}_{i-1} \circ \underbrace{s^\AP_1 \circ\ldots \circ s^\AP_1}_{n-i}\right)(A_{1,2}) = s_1^{i-1} s_2^{n-i}(A_{1,2}) \quad 1 \le i \le n.\]
Let \(a_i \defeq \Theta_n(z_i)\) for \(1 \le i \le n\).
Since finitely generated free groups are Hopfian, it suffices to show that the braids \(a_1, \ldots, a_n\) generate a free subgroup of rank \(n\).

\begin{figure}[H]
\centering
\includegraphics[width=0.3\textwidth]{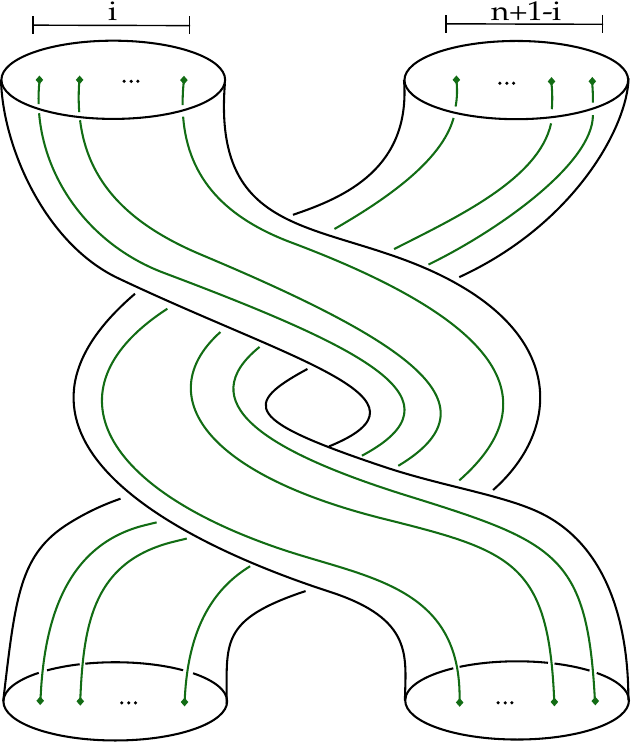}
\caption{Picture of the geometric braid \(a_i\).}
\end{figure}

\begin{Proposition}{}{imviatilde}
For every \(1 \le i \le n\), one has
\[a_i = s_1^{i-1} s_2^{n-i}(A_{1,2}) = (\widetilde A_{0,1} \ldots \widetilde A_{0,i})^{-1}.\]
\end{Proposition}
\begin{proof}
It is convenient to use explicit formulas for \(s^\AP_* = s_{*+1}\):
\begin{equation}
\label{s_formula}
s_1(A_{s,k}) = \begin{cases}
A_{0,1} A_{0,2} A_{1,2}^2, & (s,k) = (0,1), \\
A_{0,k+1}, & s = 0, k > 1, \\
A_{1,k+1}A_{2,k+1}, & s = 1, \\
A_{s+1,k+1}, & s > 1.
\end{cases} \quad s_2(A_{s,k}) = \begin{cases}
A_{0,1}, & (s,k) = (0,1), \\
A_{1,2}A_{1,3}, & (s,k) = (1,2), \\
A_{0,2} A_{0,3} A_{2,3}^2, & (s,k) = (0,2), \\
A_{s,k+1}, & s < 2 < k, \\
A_{2,k+1}A_{3,k+1}, & s = 2, \\
A_{s+1,k+1}, & s > 2.
\end{cases}
\end{equation}

First, we claim that for every \(1 \le i \le n\) one has
\[s_2^{n-i}(A_{1,2}) = A_{1,2} \ldots A_{1,n-i+2} = A_{0,1}^{-1} \in P_{n-i+2}.\]
The proof is by reverse induction on \(i\). For the base case \(i = n\), there is nothing to show. If the statement is proved for some value \(i > 1\), then
\begin{align*}
s_2^{n-i+1}(A_{1,2}) &= s_2 s_2^{n-i}(A_{1,2}) = s_2(A_{0,1}^{-1}) \quad \text{by induction} \\
&= A_{0,1}^{-1} \quad \text{by \eqref{s_formula}}.	
\end{align*}
The claim is proved.

Second, we prove the proposition.
The proof is by induction on \(i\). The base case \(i = 1\) follows from the claim (\(t_1 = 1\) and \(\widetilde A_{0,1} = A_{0,1}\)). If the statement is proved for some value \(i\) and every number of strands greater than \(i\), then
\begin{align*}
s_1^i s_2^{n-i-1}(A_{1,2}) &= s_1 s_1^{i-1} s_2^{n-i-1}(A_{1,2}) = s_1(A_{0,1} \cdot A_{0,2} t_2^2 \ldots A_{0,i} t_i^2)^{-1} \quad \text{by induction}\\
&= (A_{0,1} A_{0,2} A_{1,2}^2 \cdot A_{0,3} t_3^2 \ldots A_{0,i+1} t_{i+1}^2)^{-1} \quad \text{by \eqref{s_formula}}\\
&= (\widetilde A_{0,1} \ldots \widetilde A_{0,i+1})^{-1}.
\end{align*}
\end{proof}

Therefore,
\[\mathrm{Im}(\Theta_n) = \langle a_1, \ldots, a_n \rangle \eqcapt{\text{Prop.~\ref{t:imviatilde}}} \langle \widetilde A_{0,1}, \ldots, \widetilde A_{0,n} \rangle \eqcapt{\text{Prop.~\ref{t:bigprod}}} \langle \widetilde A_{0,2}, \ldots, \widetilde A_{0,n+1} \rangle.\]

Hence, we reduce the initial problem to the following one.
\begin{Red}{reformulation of the injectivity of \(\Theta_n\)}{}
Show that the subgroup in \(P_{n+1}\) generated by braids \(\widetilde A_{0,2}, \ldots, \widetilde A_{0, n+1}\) is free of rank \(n\).
\end{Red}

Consider the following braids \(x_1, \ldots, x_n, y_1, \ldots, y_{n-1} \in P_{n+1}\):
\begin{align*}
x_i &\defeq A_{i,n+1}, \\
y_i &\defeq A_{i,i+1} A_{i,i+2} \ldots A_{i,n}.
\end{align*}
For convenience, we set \(y_n \defeq 1\).
\begin{Lemma}{}*
One has
\[\mu_{n+1}(\widetilde A_{0,k}) = (y_{k-1} x_{k-1} y_{k-1})^{-1}, 2 \le k \le n+1.\]
\end{Lemma}
\begin{proof}
Note that for every \(1 \le k \le n+1\) one has
\begin{align*}
\mu_{n+1}(t_k) &= \mu_{n+1}(A_{1,k} \cdot A_{2,k} \ldots A_{k-1,k}) = A_{0,k-1} \cdot A_{1,k-1} \ldots A_{k-2,k-1} \\
&= (A_{k-1,k+1} \ldots A_{k-1,n+1}) = x_{k-1}^{-1} y_{k-1}^{-1}.
\end{align*}

If \(k > 1\), then
\[\mu_{n+1}(\widetilde A_{0,k}) = \mu_{n+1}(A_{0,k} \cdot t_k^2) \eqcapt{\text{Lm.~\ref{t:mucomp}}} x_{k-1} \cdot x_{k-1}^{-1} y_{k-1}^{-1} x_{k-1}^{-1} y_{k-1}^{-1} = (y_{k-1} x_{k-1} y_{k-1})^{-1}.\]
\end{proof}

Consider a subgroup \(H_n \defeq \mu_{n+1}(\mathrm{Im}(\Theta_n))\) of \(P_{n+1}\).
One has
\[H_n = \langle \mu_{n+1}(\widetilde A_{0,2}), \ldots, \mu_{n+1}(\widetilde A_{0,n+1})\rangle = \langle y_1 x_1 y_1, \ldots y_n x_n y_n \rangle.\]

Let \(G_n \defeq \langle x_1, x_2, \ldots, x_n, y_1, \ldots, y_{n-1} \rangle \subset P_{n+1}\).
The next lemma gives the presentation for the group \(G_n\). It is convenient to introduce braids \(p_1, p_2, \ldots, p_n \in P_{n+1}\) using the following formula:
\[p_i \defeq x_i x_{i+1} \ldots x_n.\]
For convenience, we set \(p_{n+1} \defeq 1\). Note that \(x_i = p_i p_{i+1}^{-1}\) for \(1 \le i \le n\).

\begin{figure}[H]
\centering
\includegraphics[width=0.75\textwidth]{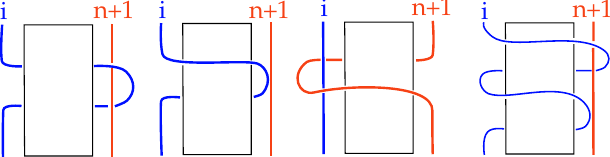}
\caption{Pictures of the geometric braids \(x_i, y_i, p_i\), and \(y_ix_iy_i\).}
\end{figure}

\begin{Lemma}{}*
The group \(G_n\) splits as a semidirect product
\[G_n = \langle x_1, \ldots, x_n \rangle \rtimes \langle y_1, \ldots, y_{n-1} \rangle \cong F_n \rtimes \ZZ^{n-1},\]
with an action of \(y_j\) on \(x_i\) given by
\[x_i^{y_j} = \begin{cases}
x_i^{x_j^{-1}}, & j < i, \\
x_i^{p_i^{-1}}, & j = i, \\
x_i, & j > i.
\end{cases}\]
\end{Lemma}
\begin{proof}
Consider a short exact sequence
\[\mathrm{Ker}(d_{n+1}) \hookrightarrow P_{n+1} \xtwoheadrightarrow{d_{n+1}} P_n.\]
It is known that \(\mathrm{Ker}(d_{n+1}) = \langle x_1, x_2, \ldots, x_n \rangle\) is free of rank \(n\), and this sequence splits with a section \(\iota \colon P_n \to P_{n+1}\) given by \(A_{i,j} \mapsto A_{i,j}\). Hence,
\begin{equation}
\label{eq:splitting}
P_{n+1} = \langle x_1, \ldots, x_n \rangle \rtimes P_n.
\end{equation}
Since \(G_n \supset \langle x_1, \ldots, x_n \rangle\), the splitting~\eqref{eq:splitting} restricts to the splitting
\[G_n = \langle x_1, \ldots, x_n \rangle \rtimes (G_n \cap \iota(P_n)) = \langle x_1, \ldots, x_n \rangle \rtimes \langle y_1, \ldots, y_{n-1} \rangle.\]

The action of \(y_j\) on \(x_i\) can be easily seen pictorially.
\begin{figure}[H]
\centering
\includegraphics[width=0.75\textwidth]{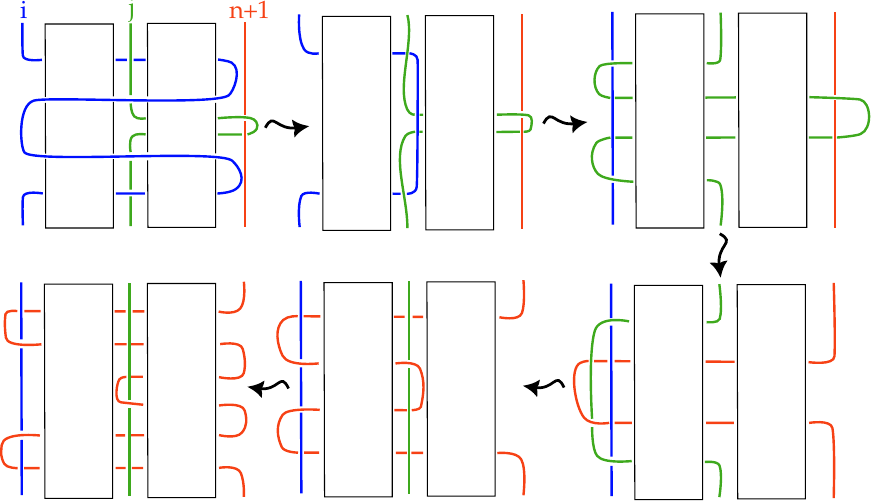}
\caption{The equality \(x_i^{y_j} = x_i^{x_j^{-1}}\) for \(j < i\).}
\end{figure}
\begin{figure}[H]
\centering
\includegraphics[width=0.6\textwidth]{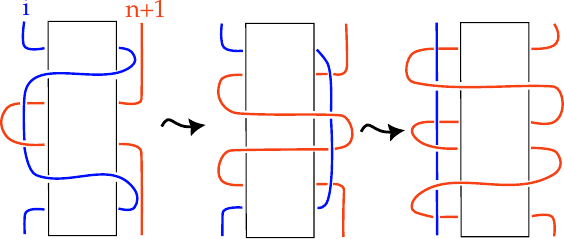}
\caption{The equality \(x_i^{y_i} = x_i^{p_i^{-1}}\).}
\end{figure}
\end{proof}

The braids \(p_1, p_2, \ldots, p_n\) form another free basis for the group \(\langle x_1, x_2, \ldots, x_n \rangle \subset P_{n+1}\).
\begin{Lemma}{}*
The action of \(y_j\) on \(p_i\) is given by
\[p_i^{y_j} = \begin{cases}
p_i^{p_{j+1} p_j^{-1}}, & j < i, \\
p_i, & j \ge i.
\end{cases}\]
\end{Lemma}
\begin{proof}
If \(j < i\), then
\begin{align*}
p_i^{y_j} = (x_i x_{i+1} \ldots x_n)^{y_j} = (x_i x_{i+1} \ldots x_n)^{x_j^{-1}} = p_i^{p_{j+1} p_j^{-1}}.
\end{align*}
If \(j = i\), then
\begin{align*}
p_i^{y_i} &= (x_i \cdot x_{i+1} \ldots x_i)^{y_i} = x_i^{(x_i x_{i+1} \ldots x_n)^{-1}} \cdot x_{i+1}^{x_i^{-1}} \ldots x_n^{x_i^{-1}}\\
 &= (p_i p_{i+1}^{-1})^{p_i^{-1}} \cdot p_{i+1}^{p_{i+1} p_i^{-1}} = p_i.
\end{align*}
If \(j > i\), then
\begin{align*}
p_i^{y_j} &= (x_i \ldots x_{j-1} \cdot x_j \cdot x_{j+1} \ldots x_n)^{y_j} = x_i \ldots x_{j-1} \cdot x_j^{(x_j x_{j+1} \ldots x_n)^{-1}} \cdot (x_{j+1} \ldots x_n)^{x_j^{-1}} \\
&= p_i p_j^{-1} \cdot (p_j p_{j+1}^{-1})^{p_j^{-1}} \cdot p_{j+1}^{p_{j+1}p_j^{-1}} = p_i.
\end{align*}
\end{proof}
Note that \(y_i x_i y_i = y_i p_i p_{i+1}^{-1} y_i\) for every \(1 \le i \le n\).

\begin{Red}{reformulation of the injectivity of \(\Theta_n\)}{problem2}
Let \(G_n\) be a group given by the following presentation:
\begin{align*}
G_n = \langle p_1, \ldots, p_n, y_1, \ldots, y_{n-1} \mid [y_i, y_j] = 1; \; [p_i, y_j] = 1, \; j \ge i; \; [p_i, y_j p_j p_{j+1}^{-1}] = 1, \; j < i \rangle.
\end{align*}
Also, let \(y_n \defeq 1, p_{n+1} \defeq 1\).
Prove that the subgroup \(H_n = \langle y_1 p_1 p_2^{-1} y_1, \ldots, y_n p_n p_{n+1}^{-1} y_n\rangle\) of \(G_n\) is free of rank \(n\).
\end{Red}

\section{Solution to Problem~\ref{t:problem2}}
Let \(K_n\) be a group given by the following presentation:
\[K_n = \langle p_1, \ldots, p_{n-1}, y_1, \ldots, y_{n-1} \mid [y_i, y_j] = 1; \; [p_i, y_j] = 1, \; j \ge i; \; [p_i, y_j p_j p_{j+1}^{-1}] = 1, \; j < i \rangle.\]
As before, we set \(x_i \defeq p_i p_{i+1}^{-1}\) for \(1 \le i < n-1\).
Consider a subgroup
\[A_n \defeq \langle y_1x_1, \ldots, y_{n-2}x_{n-2}, y_{n-1} p_{n-1} \rangle \subset K_n.\]
Then \(G_n\) can be obtained from \(K_n\) as an HNN extension
\[G_n \cong \langle K_n, p_n \mid a^{p_n} = a, \; a \in A_n \rangle.\]
We set \(x_{n-1} \defeq p_{n-1} p_n^{-1}, x_n \defeq p_n \in G_n\).
\begin{Lemma}{}*
One has
\begin{align*}
H_n &\eqbydef \langle y_1 x_1 y_1, \ldots, y_{n-2} x_{n-2} y_{n-2}, y_{n-1} x_{n-1} y_{n-1}, x_n\rangle \\
&= \langle y_1 x_1 y_1, \ldots, y_{n-2} x_{n-2} y_{n-2}, y_{n-1} p_{n-1} y_{n-1}, p_n \rangle.
\end{align*}
\end{Lemma}
\begin{proof}
Note that
\[p_n^{y_{n-1} p_{n-1}} = p_n^{p_n p_{n-1}^{-1} p_{n-1}} = p_n.\]
Hence,
\[g \defeq y_{n-1} x_{n-1} y_{n-1} = y_{n-1} p_{n-1} p_n^{-1} y_{n-1} = p_n^{-1} (y_{n-1} p_{n-1} y_{n-1}).\]
Therefore, we can replace \(g\) with \(p_n g = y_{n-1} p_{n-1} y_{n-1}\) in the initial set of generators.
\end{proof}
The advantage of the new generating set is that all but one of its members are in \(K_n\). Let \(L_n\) be a subgroup of \(G_n\) given by
\[L_n = \langle y_1 x_1 y_1, \ldots, y_{n-2} x_{n-2} y_{n-2}, y_{n-1} p_{n-1} y_{n-1} \rangle.\]
Then,
\[L_n \subset K_n, \; H_n = \langle L_n, p_n \rangle.\]
We aim to prove that \(H_n\) is free of rank \(n\) by induction. Since \(d_n(L_n) = H_{n-1}\), the group \(L_n\) would be free of rank \(n-1\) by the inductive hypothesis. Hence, it would be enough to prove that the natural map
\[L_n * \langle p_n \rangle \to G_n = \langle K_n, p_n \mid [p_n, a] = 1, \; a \in A_n \rangle\]
is injective.

The following lemma is an immediate corollary of the Britton's normal form \cite[\S 4]{LS01}.
\begin{Lemma}{}{britton}
Let \(K\) be a group, and let \(\varphi \colon A \to B\) be an isomorphism between some subgroups \(A, B \subset G\).
Let
\[G = \langle K, t \mid a^t = \varphi(a),\; a \in A \rangle\]
be an HNN extension. If a subgroup \(L \subset K\) is such that \(L \cap A = L \cap B = 1\), then the natural map \(L * \langle t \rangle \to G\) is injective. \qed
\end{Lemma}

\begin{Theorem}{}{main}
For every \(k \ge 1\) one has
\begin{enumerate}[label=(\arabic*)]
\item \(L_k \cap A_k = 1\);
\item \(H_k\) is free of rank \(k\);
\item \(H_k \cap \langle A_k, p_k \rangle = \langle p_k \rangle\).
\end{enumerate}
\end{Theorem}
\begin{proof}
We will prove this theorem by induction on \(k\). The scheme of the proof resembles a long-exact sequence:
\begin{center}
\tikzset{every picture/.style={line width=0.75pt}}
\begin{tikzpicture}[x=0.75pt,y=0.75pt,yscale=-1,xscale=1]
\draw   (372.71,105) -- (372.71,91.83) .. controls (372.71,86.56) and (376.99,82.29) .. (382.25,82.29) -- (605.9,82.29) .. controls (611.17,82.29) and (615.44,86.56) .. (615.44,91.83) -- (615.44,98.19) -- (617.71,98.19) -- (612.04,105) -- (606.36,98.19) -- (608.63,98.19) -- (608.63,91.83) .. controls (608.63,90.32) and (607.41,89.1) .. (605.9,89.1) -- (382.25,89.1) .. controls (380.75,89.1) and (379.53,90.32) .. (379.53,91.83) -- (379.53,105) -- cycle ;
\draw   (489.71,126) .. controls (489.71,121.33) and (487.38,119) .. (482.71,119) -- (385.66,119) .. controls (378.99,119) and (375.66,116.67) .. (375.66,112) .. controls (375.66,116.67) and (372.33,119) .. (365.66,119)(368.66,119) -- (272.71,119) .. controls (268.04,119) and (265.71,121.33) .. (265.71,126) ;
\draw   (99.71,151) .. controls (99.71,155.67) and (102.04,158) .. (106.71,158) -- (221.65,158) .. controls (228.32,158) and (231.65,160.33) .. (231.65,165) .. controls (231.65,160.33) and (234.98,158) .. (241.65,158)(238.65,158) -- (351.71,158) .. controls (356.38,158) and (358.71,155.67) .. (358.71,151) ;
\draw   (230.71,170) -- (230.71,184.68) .. controls (230.71,188.16) and (233.53,190.97) .. (237.01,190.97) -- (444.32,190.97) .. controls (447.8,190.97) and (450.62,188.16) .. (450.62,184.68) -- (450.62,180.49) -- (452.71,180.49) -- (447.47,174.19) -- (442.23,180.49) -- (444.32,180.49) -- (444.32,184.68) .. controls (444.32,184.68) and (444.32,184.68) .. (444.32,184.68) -- (237.01,184.68) .. controls (237.01,184.68) and (237.01,184.68) .. (237.01,184.68) -- (237.01,170) -- cycle ;
\draw   (205.71,137.79) -- (243.49,137.79) -- (243.49,135.29) -- (250,140.29) -- (243.49,145.29) -- (243.49,142.79) -- (205.71,142.79) -- cycle ;
\draw (99,132) node [anchor=north west][inner sep=0.75pt]   [align=left] {(1) for $\displaystyle k=n$};
\draw (262,130) node [anchor=north west][inner sep=0.75pt]   [align=left] {(2) for $\displaystyle k=n$};
\draw (402,130) node [anchor=north west][inner sep=0.75pt]   [align=left] {(3) for $\displaystyle k=n$};
\draw (542,130) node [anchor=north west][inner sep=0.75pt]   [align=left] {(1) for $\displaystyle k=n+1$};
\draw (205,116) node [anchor=north west][inner sep=0.75pt]   [align=left] {\textbf{Step I}};
\draw (325,194) node [anchor=north west][inner sep=0.75pt]   [align=left] {\textbf{Step II}};
\draw (473,61) node [anchor=north west][inner sep=0.75pt]   [align=left] {\textbf{Step III}};
\end{tikzpicture}
\end{center}

\textbf{Step I}. Obviously, (1) implies (2) by Lemma~\ref{t:britton}.

\textbf{Step II}. Assume (1) and (2) for \(k = n\). Then
\[H_n \cap \langle A_n, p_n \rangle = \langle L_n, p_n \rangle \cap \langle A_n, p_n \rangle.\]
Since \(p_n\) acts trivially by conjugation on \(\langle A_n, p_n \rangle\), every element \(x\) from the intersection must commute with \(p_n\). On the other hand, since \(p_n\) is an element of the free basis of \(H_n\), \(x\) must be some power of \(p_n\). This implies (3).

\textbf{Step III}. Assume (2) and (3) for \(k = n\). Consider a homomorphism \(d_{n+1} \colon P_{n+2} \to P_{n+1}\). It is easy to see that \(d_{n+1}(A_{n+1}) = \langle A_n, p_n \rangle\), and \(d_{n+1}(L_{n+1}) = H_n\). Also note that the restriction of \(d_{n+1}\) to \(L_{n+1}\) must be injective, since its image, \(H_n\), is free of rank \(n\) by (2). Hence,
\begin{align*}
L_{n+1} &\cap A_{n+1} \subset d_{n+1}^{-1}(H_n \cap \langle A_n, p_n \rangle) \cap L_{n+1} \quad \text{by (1)} \\
&= d_{n+1}^{-1}(\langle p_n \rangle) \cap L_{n+1} \quad \text{by inductive hypothesis} \\
&= \langle y_n p_n y_n \rangle \quad \text{since}\; d_{n+1}|_{L_{n+1}} \;\text{is mono}.
\end{align*}
It is left to show that \(A_{n+1}\) and \(\langle y_n p_n y_n \rangle\) intersect trivially.
Consider a homomorphism
\[f = d_1 \circ d_2 \circ \ldots \circ d_{n-1} \colon P_{n+2} \to P_3.\]
The subgroup \(\langle y_n p_n y_n \rangle\) maps injectively to \(\langle y_1x_1x_2y_1\rangle\) under \(f\). Also, one has \(f(y_i x_i) = 1\) whenever \(1 \le i < n\), and \(f(y_nx_nx_{n+1}) = y_1x_1x_2\). Hence, \(f(A_{n+1}) \subset \langle y_1x_1x_2 \rangle\). Therefore,
\[\langle y_n p_n y_n \rangle \cap A_{n+1} \xhookrightarrow{f} \langle y_1 x_1 x_2 y_1 \rangle \cap \langle y_1x_1x_2 \rangle.\]
The result follows since the subgroups
\begin{equation*}
\langle y_1 x_1 x_2 y_1 = A_{1,2} A_{1,3} A_{2,3} A_{1,2} \rangle \; \text{and} \; \langle y_1x_1x_2 = A_{1,2} A_{1,3} A_{2,3} \rangle  
\end{equation*}
do not intersect in \(P_3\) (for example, since \(Z(P_3) = \langle A_{1,2}A_{1,3} A_{2,3}\rangle\), the latter subgroup is central, and the former is not).
\end{proof}
The part (2) of Theorem~\ref{t:main} provides a solution to Problem~\ref{t:problem2}, and thereby proves Theorem~\ref{t:cohen-wu-embedding}.

\end{document}